\input amssym.def
\input amssym
\magnification=1200
\parindent0pt
\hsize=16 true cm
\baselineskip=13  pt plus .2pt
$ $

\def\Z{\Bbb Z}

\centerline {\bf On the determination of knots by their cyclic unbranched
coverings}

\bigskip

\centerline {Bruno P. Zimmermann}

\medskip

\centerline {Universit\`a degli Studi di Trieste}
\centerline {Dipartimento di Matematica e Informatica}
\centerline {34100 Trieste, Italy}
\centerline {zimmer@units.it}

\bigskip \bigskip

\vskip 1cm

Abstract. {\sl We show that, for any  prime $p$, a knot $K$ in $S^3$ is
determined by its $p$-fold cyclic unbranched covering. We also investigate when
the $m$-fold cyclic unbranched covering of a knot in $S^3$ coincides with the
$n$-fold cyclic unbranched covering of another knot, for different coprime
integers $m$ and $n$.}

\bigskip \bigskip

{\bf 1. Introduction}

\medskip

There is an extensive literature on the determination of knots in $S^3$ by their
$p$-fold cyclic {\it branched} coverings; for example, the case of odd prime
numbers $p$ is considered in [Z1] for hyperbolic knots and in [BP] for arbitrary
prime knots, the case of 2-fold branched coverings is considered in [Rn],[MZ] and
[K] (see also the survey [Z2]). On the other hand, less seems to be known for the
case of cyclic {\it unbranched} coverings (that is, of the complements of the
knots). For a knot $K$ in $S^3$, we denote by $M_p(K)$ the
$p$-fold cyclic unbranched covering of its complement $M_1(K) = S^3 - N(K)$
($S^3$ minus the interior of a regular neighbourhood of the knot), so $M_p(K)$ is
a compact orientable 3-manifold with a torus boundary.  The basic case here is,
of course, the case $p=1$ or the fact that a knot $K$ in $S^3$ is determined by
its complement ([GoL]). In the present paper, we study the case of
primes $p>1$ and prove the following

\bigskip

{\bf Theorem 1} {\sl For any prime $p$, a knot
$K$ in $S^3$ is determined by its  $p$-fold cyclic unbranched covering $M_p(K)$
(i.e., any other knot with the same  $p$-fold cyclic unbranched covering is
equivalent to  $K$).}

\bigskip

Here two unoriented knots are equivalent if there is a diffeomorphism of $S^3$
which maps one to the other.  Concerning different branching orders, we introduce
the following

\bigskip

{\bf Abelian Construction}.  Let $m$ and $n$ be distinct positive integers. Let
$\bar K$ be a knot in a lens space $L$  which represents a generator of
$\pi_1(L) \cong \Bbb Z_n$, and denote by $K$ the knot which is the
preimage of $\bar K$ in the universal covering $S^3$ of $L$. Similarly, let $K'$
be the preimage of a knot $\bar K'$ in a lens space $L'$, representing a
generator of $\pi_1(L') \cong \Bbb Z_m$, in  the universal covering $S^3$ of
$L'$.  Suppose that
$\bar K$ and $\bar K'$ have homeomorphic complements $L- N(\bar K) = L'-
N(\bar K')$; then the $m$-fold cyclic unbranched covering of $K$ coincides with
the $n$-fold cyclic unbranched covering of $K'$, i.e.  $M_m(K) = M_n(K')$  (both
are equal to a regular unbranched ($\Bbb Z_m \times \Bbb Z_n$)-covering of $L-
N(\bar K) = L'- N(\bar K')$).

\bigskip

Concentrating mainly on the basic case of hyperbolic knots, the following
holds.

\bigskip

{\bf Theorem 2} {\sl i) Let $K$ be a hyperbolic knot and $K'$ be any knot in $S^3$
such that  $M_m(K) = M_n(K')$, for coprime positive integers $m$ and $n$. Then
$K$ and $K'$ are obtained by the Abelian Construction. The same remains true for
arbitray knots $K$ and $K'$ if $m$ and $n$ are different prime numbers.

iii) Let $M$ be a compact orientable 3-manifold whose boundary is a torus and
whose interior has a complete hyperbolic structure of finite volume. There are at
most three coprime positive integers $m$ such that $M$ is the $m$-fold cyclic
unbranched covering of a knot $K$ in $S^3$.}

\bigskip

We think that part i) of Theorem 2 remains true for arbitrary knots $K$ and
$K'$ and coprime integers $m$ and $n$.  We note that the present formulation of
part i) of Theorem 2 uses the recent geometrization of  free cyclic group
actions on $S^3$ after Perelman ([P1], [P2]; without this, one concludes that,
in the Abelian Construction, $L$ and $L'$ are 3-manifolds with finite cyclic
fundamental groups whose universal covering is  $S^3$.

\medskip

Let $K_0$ be a knot in $S^3$ which admits two non-trivial lens space surgeries
$L$ and $L'$, with $\pi_1(L) \cong \Bbb Z_n$ and $\pi_1(L') \cong \Bbb Z_m$. The
cores of the surgered solid tori give two knots $\bar K$ and $\bar K'$ in the
lens spaces $L$ and $L'$ whose complements coincide with the complement of
$K_0$ in $S^3$.  By the Abelian Construction one obtains two knots $K$ and $K'$
in $S^3$ such that  $M_m(K) = M_n(K') = M_{mn}(K_0)$. Note also that  $M_1(K) =
M_n(K_0)$ and $M_1(K') =  M_m(K_0)$. Examples of hyperbolic knots in  $S^3$ with
two non-trivial lens space surgeries can be found in [Be], [FS].

\medskip

The case of different branching orders for cyclic {\it branched} coverings has
been considered in [RZ],[Z2] for hyperbolic knots, see also [BPZ]; one main
difference is that the isometry group of such a cyclic branched covering (a
closed hyperbolic 3-manifold) may, in principle, be much more complicated than
that of a cyclic unbranched covering (a 3-manifold with torus-boundary). In
particular, the situation for cyclic branched coverings is not well understood in
the case of non-solvable isometry groups; for example, it is not clear which
finite non-abelian simple groups may occur as groups of isometries of a cyclic
branched covering.

\medskip

We thank  S. Friedl, L. Paoluzzi and M. Scharlemann for helpful correspondence,
and A. Reid for bringing to our attention the references [GHH] and [ReW] on the
related topic of commensurability classes of hyperbolic knot complements.

\bigskip

{\bf 2. Proof of Theorem 1}

\medskip

We start with the following

\bigskip

{\bf Lemma 1} {\sl  For a prime $p$, let $G = {\Bbb Z}_p \times {\Bbb Z}_p$  be a
finite group of orientation-preserving diffeomorphisms of a mod $p$ homology
3-sphere
$M$ (i.e., for homology with coefficients in the integers mod $p$). Then either
there are exactly two subgroups $\Bbb Z_p$ of
$G$ with non-empty fixed point set (two disjoint circles) or, if $p = 2$,  all
three involutions in $G$ may have nonempty fixed point set  (three circles
intersecting in exactly two points).}

\medskip

{\sl Proof.}  By  Smith fixed point theory, $G$ does not act freely, and each
element of $G$ has empty or connected fixed point set (see [Br, Theorems 7.9 e
8.1]). Let $X$ be a nontrivial cyclic subgroup of
$G$ with nonempty fixed point set $K$ which is a circle.  Because $G$ is abelian,
$K$ is invariant under the action of $G$.  The projection $\bar G$ of
$G$ to $\bar M := M/X$ is a cyclic group leaving invariant the projection $\bar
K$ of $K$. It is easy to see that also $\bar M$ is a  mod $p$ homology 3-sphere,
and hence $\bar G$ has empty or connected fixed point set.

\medskip

Suppose that there is another nontrivial cyclic subgroup $X' \ne X$ of $G$ with
nonempty fixed point set $K'$ different from $K$.  If $K$ and $K'$ intersect then
they intersect in exactly two points and $p=2$ (because $K$ is invariant under
$X'$), and consequently we are in the second case of the Lemma. Hence we can
assume that $K$ and $K'$ do not intersect. Note that
$K'$ is invariant under $X$ which acts as a group of rotations on $K'$. As the
fixed point set of $\bar G$ is connected (a circle) it consists of the projection
of $K'$. The preimage of this projection is exactly $K'$ which implies that $X$
and $X'$ are the only nontrivial cyclic subgroups of $G$ with non-empty fixed
point set. Thus we are in the first case of the Lemma.

\medskip

Now suppose that $X$ is the only cyclic subgroup of $G$ with nonempty fixed point
set.  Then $\bar G$ acts freely on the mod $p$ homology 3-sphere $\bar M$. Let $N
:= \bar M/\bar G$ be the quotient and $L$ the projection of $\bar K$ to $N$. Note
that $M - K$ is a regular unbranched  covering of $N - L$, with covering group $G
\cong {\Bbb Z}_p \times {\Bbb Z}_p$. By [H, p.92], $H_1(N - L;\Z_p)$ is
isomorphic to  $\Bbb Z_p$; this implies that $N - L$ has no abelian covering with
covering group ${\Bbb Z}_p \times {\Bbb Z}_p$ which is a contradiction

\medskip

This finishes the proof of Lemma 1.

\bigskip

Starting with the {\it Proof of Theorem 1} now, the $p$-fold cyclic unbranched
covering $M = M_p(K)$ of the knot $K$ is a compact 3-manifold whose boundary is a
torus; we denote by $C \cong \Bbb Z_p$ the cyclic covering group acting
freely on $M$. Suppose that $M$ is also the
$p$-fold cyclic unbranched covering of another knot $K'$, i.e. $M = M_p(K')$,
with covering group $C' \cong \Bbb Z_p$.  We denote by $B_p(K)$ the $p$-fold
cyclic {\it branched} covering of $K$, so $M = B_p(K) - N(\tilde K)$ where
$N(\tilde K)$ denotes the interior of a regular neighbourhood of the preimage
$\tilde K$ of $K$ in $B_p(K)$. The action of the covering group $C$ on $M$ extends
to an action of $C$ on $B_p(K)$ with fixed point set $\tilde K$, giving the
covering group of the
$p$-fold cyclic branched covering of $K$ which will be denoted also by $C$.
Similarly, $C'$ extends to the $p$-fold cyclic branched covering $B_p(K')$ of
$K'$ fixing the preimage $\tilde K'$ of $K'$.

\medskip

We will show that, up to conjugation, the covering groups $C$
and $C'$ commute; Theorem 1 then follows from the following

\bigskip

{\bf Lemma 2} {\sl  Suppose that the covering groups $C$ and $C'$ commute,
generating a group $G = C \oplus C' \cong \Bbb Z_p \times \Bbb Z_p$ of
diffeomorphisms of $M$. Then the knots $K$ and $K'$ are equivalent.}

\medskip

{\it Proof.}   The action of $C'$ on
$M$ extends to a free action of $C'$ on the $p$-fold cyclic branched covering
$B_p(K)$ of $K$ (unless $C=C'$: in this case, $K$ and $K'$ have homeomorphic
complements and hence are equivalent).  We can assume that the actions of
$C$ and $C'$ commute also on $B_p(K)$ and hence generate a group $G = C \oplus C'
\cong  \Bbb Z_p \times \Bbb Z_p$ of diffeomorphisms of $B_p(K)$.

\medskip

The $p$-fold cyclic branched covering $B_p(K)$ of a knot $K$ in $S^3$ is a mod
$p$ homology 3-sphere (see e.g. [Go]). By Lemma 1, there are exactly two
subgroups $\Z_p$ of $G$ with non-empty (connected) fixed point set; one of these
is  $C$ which fixes $\tilde K$, and we denote by $\tilde A$ the fixed point set
of the other which is contained in $M$. Now $G$ projects to a group
$H \cong G/C  \cong \Bbb Z_p$ of symmetries of $(B_p(K), \tilde K)/C = (S^3, K)$,
with non-empty fixed point set $A$ disjoint from $K$ (the projection of
$\tilde A$). Hence $K$ has cyclic period $p$.

\medskip

By the positive solution of the Smith conjecture, $H$ acts by standard rotations
on the 3-sphere, so $S^3/H$ is again the 3-sphere. The group $H$ acts by
rotations along $K$ and maps a meridian of $K$ to $p$ disjoint meridians of
$K$. Its fixed point set is the projection $\bar A$ of $\tilde A$ resp. $A$, so
$(S^3,K)$ is the cyclic branched covering of $(S^3,\bar K)$ branched along the
trivial knot $\bar A$ in $S^3$.

\medskip

Now $(S^3 - N(K))/H = S^3 - N(\bar K) = M/G$, and by symmetry and with analogous
notation,  $(S^3 - N(K'))/H' = S^3 - N(\bar K') = M/G$, so
$S^3 - N(\bar K) = S^3 - N(\bar K')$ (where $H' \cong G/C' \cong \Bbb Z_p$).
Hence $(S^3,\bar K')$ is obtained from $(S^3,\bar K)$ by $1/n$-surgery on
$\bar K$, for some integer $n$, and this surgery transforms also
$\bar A$ into $\bar A'$ (note that the projections  $\bar A$ and $\bar A'$ of
$\tilde A \subset M$ coincide as subsets of  $M/G = S^3 - N(\bar K) = S^3 -
N(\bar K')$).

\medskip

If the surgery on $\bar K$ is trivial (i.e., $n=0$) then $\bar K' = \bar K$ and
$\bar A' =  \bar A$, so $K' = K$ and Lemma 2 is proved in this case.

\medskip

In the case of non-trivial surgery, since the result of the surgery is again
the 3-sphere, $\bar K$ has to be a trivial knot by [GoL], and hence $\bar K
\cup \bar A$ is a hyperbolic link of two unknotted components.  By [M,
Corollary 3], a non-trivial $1/n$-surgery on one component $\bar K$ of the
hyperbolic link $\bar K \cup \bar A$ transforms the other component $\bar A$
into a non-trivial knot (obtained from $\bar A$ by twisting $n$ times around a
spanning disk for $\bar K$), and hence $\bar A'$ is non-trivial. This
contradicts the Smith conjecture since the $p$-fold cyclic covering of $S^3$
branched along  $\bar A'$ is the 3-sphere.

\medskip

This finishes the proof of Lemma 2.

\bigskip

The proof of Theorem 1 follows now from the following

\bigskip

{\bf Lemma 3} {\sl Up to conjugation, the covering groups $C$ and $C'$ commute.}

\medskip

{\it Proof.}  Since the proof is much easier if $K$ is a hyperbolic knot we
will give the proof first for this case.

\bigskip

i) Suppose that $K$ is hyperbolic. Then $M$ is a hyperbolic 3-manifold of finite
volume, with one torus-boundary or cusp. It
is a consequence of Mostow rigidity and Waldhausen's theorem for Haken
3-manifolds ([W]) that also
$K'$ is hyperbolic, and hence we can assume that both cyclic covering groups $C$
of $K$ and $C'$ of  $K'$ act by isometries on $M$.  The covering groups $C$ and
$C'$ act freely on $M$, and hence by euclidean rotations on the boundary torus
of the hyperbolic 3-manifold $M$ (corresponding to a cusp of the hyperbolic
3-manifold $M$; the rotations lift to translations of the euclidean horospheres
corresponding to the boundary torus of $M$ in the universal covering of $M$). It
follows that the groups $C$ and $C'$ of isometries of $M$ commute elementwise
(because they commute on the boundary torus of $M$).

\medskip

This finishes the proof of Lemma 3 and Theorem 1 in the case where $K$ is
hyperbolic.

\bigskip

ii)  Now let $K$ be an arbitrary knot. We will apply the methods in [BP] to show
that $C$ and $C'$ commute, up to conjugation.  We consider the JSJ-
or torus-decomposition of $M$ and the graph $\Gamma$ dual to this
decomposition; note that $\Gamma$ is a tree since $B_p(K)$ is a mod $p$ homology
sphere. By the equivariant torus-decomposition, we can assume that both
$C$ and $C'$ respect the decomposition and  are geometric on the pieces of the
decomposition (i.e., isometries on the hyperbolic pieces, fiber-preserving on the
Seifert fibered pieces; these pieces correspond to the vertices of the graph
$\Gamma$, the decomposing tori to the edges).  Then $C$ and $C'$ induce a finite
group $G$ of automorphisms of the tree $\Gamma$ which fixes the vertex
corresponding to the piece containing the boundary torus of $M$. Let
$\Gamma_f$ denote the subtree of $\Gamma$ whose vertices and edges are fixed by
every element of $G$. Now one shows as in [BP] that, up to conjugation by
diffeomorphisms of $M$, one can assume that $C$ and $C'$ commute elementwise on
the submanifold $M_f$ of $M$ corresponding to the subtree $\Gamma_f$ of
$\Gamma$  (in [BP], only the case of odd primes $p$ is considered; however, in
our situation the methods work equally well for $p=2$).

\medskip

Denote by $M_c$  the maximal connected submanifold of $M$ corresponding to a
subtree $\Gamma_c$ of $\Gamma$ containing $\Gamma_f$ on which $C$ and $C'$
commute, up to conjugation by diffeomorphisms of $M$. If $M_c=M$ we are done, so
suppose that $M_c \ne M$. Consider a boundary torus $T$ of $M_c$ connecting $M_c$
with a piece $U$ of the decomposition of $M-M_c$. By the proof of [BP, Claim 9],
we can assume that the orbit of $T$ under both $C$ and $C'$ consists of the same
$p$ tori (in all other situations, commutativity of $C$ and $C'$ can be extended
to the $G$-orbit of $U$, contradicting maximality of $M_c$). The torus $T$
projects to a torus $\bar T$  of the torus-decomposition of the complement of $K$
which seperates $S^3$ into a solid torus (containing $K$) and a knot space (the
complement of a knot, containing the projection of $U$); in particular, there is a
well-defined meridian-longitude system on $\bar T$, and also on each torus of the
$G$-orbit of $T$ which is invariant under the actions of $C$ and $C'$. Now one
replaces the knot complement by a solid torus such that one obtains again the
3-sphere, and similarly performs $C$- and $C'$-equivariant surgery on $T$ and its
images under $C$ and $C'$. Moreover, if the piece of $M_c$ containing $T$ is
hyperbolic, one does the surgery such that the resulting 3-manifold is still
hyperbolic and the central curve of the added solid torus is a shortest geodesic;
if it is Seifert fibered, one creates an exceptional fiber of high order by the
surgery.  Doing this for each such boundary torus $T$ of $M_c$, one obtains a
3-manifold $\hat M_c$ with induced actions of $C$ and $C'$ which are still the
covering groups of two $p$-fold cyclic unbranched coverings of two knots in
$S^3$; moreover, by construction $C$ and $C'$ commute on $\hat M_c$.  By Lemma 2,
the two corresponding knots are equivalent, hence the actions of $C$ and $C'$ on
$\hat M_c$ are conjugate; by the choice of the surgeries, the conjugating
diffeomorphism restricts to $M_c$ and then extends to $M$, hence we can assume
that actions of $C$ and $C'$ coincide on $M_c$.  Now by [BP, Lemma 10] the
actions of $C$ and $C'$ coincide on $M$, up to conjugation,  hence $K$ and $K'$
have homeomorphic complements and are equivalent.

\medskip

This finishes the proof of Lemma 3 and of Theorem 1 in the general case.

\bigskip

{\bf 2. Proof of Theorem 2}

\medskip

The proof is along similar lines. For the proof of part i) of Theorem 2, let
$M = M_m(K) = M_n(K')$ and denote by $C\cong \Bbb Z_m$ and $C'\cong \Bbb Z_n$
the two covering groups. As in case i) of the proof of Lemma 3, we can assume
that the covering groups act by hyperbolic isometries, commute and generate a
group $G=C \oplus C' \cong \Bbb Z_m \oplus \Bbb Z_n \cong \Bbb Z_{mn}$ of
isometries of
$M$. Both groups $C$ and $C'$ extend to the $m$-fold cyclic branched covering
$B_m(K)$ of $K$, and $C$ fixes pointwise the preimage $\tilde K$ of $K$.
The group $G$ act freely on $M$, and the only non-trivial subgroup of $G$ with
non-empty fixed point set in $B_m(K)$ is $C$. Then $G$ projects to a
cyclic group $H \cong G/C  \cong \Bbb Z_n$ acting freely on $(B_m(K),
\tilde K)/C = (S^3, K)$, so $H$ is a group of free symmetries of $K$ of order
$n$.

\medskip

Similarly, $G$ projects to a group $H' \cong G/C' \cong \Bbb Z_m$  acting
freely on $(B_n(K'), \tilde K')/C = (S^3, K')$. By the geometrization of
free cyclic group actions on $S^3$, the quotients $S^3/H$ and $S^3/H'$ are lens
spaces.  The quotients $(S^3,K)/H = (L,\bar K)$ and $(S^3,K')/H' = (L',\bar
K')$ define knots $\bar K$ and $\bar K'$ in lens spaces $L$ and $L'$ such that
$L - N(\bar K) = L' - N(\bar K') = M/G$, hence $K$ and $K'$ are obtained by the
Abelian Construction.

\medskip

This finishes the proof of Theorem 2i) in the case where $K$ is hyperbolic.
Now let $K$ be an arbitrary knot and assume that $m$ and $n$ are different
primes. We want to show that again the covering groups of the two knots commute,
up to conjugation; then the proof finishes as in the hyperbolic case. Now, in the
situation of different primes $m$ and $n$, the methods in [BPZ] apply and as in
the proof of case ii) of Theorem 1  one obtains commutativity on the
submanifold $M_f$ corresponding to the subtree $\Gamma_f$ of
$\Gamma$ fixed by the two covering groups, and then one pushes out commutativity
to all of $M$ (this is, in fact, easier than in the situation of two equal primes
considered in [BP] (as applied in the proof of Theorem 1) where in principle some
obstruction against commutativity may arise).  We remark that these methods
probably can be generalized to prove Theorem 2i) for the case of arbitrary knots
and arbitrary coprime integers.

\bigskip

For the proof of part ii) of the theorem, suppose that $M$ is the $n_i$-fold
cyclic unbranched covering of knots $K_i$ in $S^3$, for pairwise coprime
positive integers $n_1, \ldots, n_\alpha$. Denoting by $n = n_1 \ldots
n_\alpha$ their product, the manifold $M$ has now a free action of $G \cong
\Bbb Z_n$. This $G$-action on $M$ induces a free action of $G_i \cong \Bbb
Z_{n/n_i}$ on  $S^3 - N(K_i)$ which extends to a free action on $S^3$. The
quotient $S^3/G_i$ is a lens space $L_i$, with fundamental group  $\pi_1(L_i)
\cong G_i$, which contains the projection $\bar K_i$ of $K_i$. Now $L_i -
N(\bar K_i) = M/G$, so all lens spaces $L_i$ are obtained by surgery (Dehn
filling) on $M/G$. By [CGLS], there are at most three surgeries on a compact
hyperbolic 3-manifold of finite volume and with a single torus-boundary resulting
in lens spaces, hence $\alpha \le 3$.

\medskip

This finishes the proof of Theorem 2.

\bigskip  \bigskip

\centerline {\bf References}

\bigskip

\item {[Be]}  J. Berge, {\it  The knots in $D^2 \times S^1$ which have nontrivial
Dehn surgeries that yield $D^2 \times S^1$.} Topology Appl. 38,  1-19 (1991)

\item {[BP]}  M.Boileau, L.Paoluzzi, {\it  On cyclic branched coverings of prime
knots.}    To appear in  J. Topol.

\item {[BPZ]}  M.Boileau, L.Paoluzzi, B.Zimmermann, {\it  A characterization of
$\Bbb S^3$ among homology sphe-res.}  Geometry and Topology Monographs 14,
83-103  (2008)  (arXiv:math.GT/0606220)

Dedicated to the memory of Heiner Zieschang

\item {[Br]} G.Bredon, {\it Introduction to Compact Transformation Groups.}
Academic Press, New York 1972

\item {[CGLS]} M.Culler, C.McA.Gordon, J.Luecke, P.B.Shalen, {\it Dehn surgery on
knots.}  Ann. Math. 125, 237-300 (1987)

\item {[FS]} R.Fintushel, R.J.Stern,  {\it Constructing lens spaces by surgery on
knots.}  Math. Z. 175, 33-51  (1980)

\item {[GHH]} O. Goodman, D. Heard, C. Hodgson,  {\it Commensurators of cusped
hyperbolic manifolds.} See http://www.ms.unimelb.edu.au  for an electronic version

\item {[Go]} C.McA.Gordon, {\it Some aspects of classical knot theory.} Knot
Theory.  Proceedings, Plans-sur-Bex, Switzerland (J.C.Hausmann, ed.). Lect.Notes
Math. 685, 1-60,  Springer 1977

\item {[GoL]} C.McA.Gordon, J.Luecke,  {\it Knots are determined by their
complement.}  J. Amer. Math. Soc. 2, 371-415 (1989)

\item {[H]} R.Hartley, {\it Knots with free period.}  Can. J. Math. 33, 91-102
(1981)

\item {[K]} A. Kawauchi, {\it  Topological imitations and
Reni-Mecchia-Zimmermann's conjecture}. Ky-ungpook Math. J. 46, 1-9  (2006)

\item {[M]} Y. Mathieu, {\it Unknotting, knotting by twists on disks and property
(P) for knots in $S^3$.}  Knots 90. Proc. Int. Conf. Knot Theory and Related
Topics, Osaka 1990  (A. Kawauchi, ed.), 93-102, de Gruyter 1992

\item {[MZ]} M. Mecchia, B. Zimmermann,  {\it The number of knots and links with
the same 2-fold  branched covering.}  Quart. J. Math. 55, 69-76 (2004)

\item {[P1]}  G. Perelman,  {\it  The entropy formula for the Ricci flow and its
geometric applications.}  arXiv:math.DG/0211159

\item {[P2]}  G. Perelman,  {\it Ricci flow with surgery on three-manifolds.}
arXiv:math.DG/0303109

\item {[ReW]} A.W. Reid, G.S. Walsh  {\it Commensurability classes of 2-bridge
knot complements.}  arXiv: math.GT/0612473

\item {[Rn]} M.Reni, {\it On $\pi$-hyperbolic knots with the same 2-fold branched
coverings.}  Math. Ann. 316, 681-697  (2000)

\item {[RZ]} M.Reni, B.Zimmermann, {\sl Hyperbolic 3-manifolds  as cyclic
branched coverings}.  Comment. Math. Helv. 76,  300-313 (2001)

\item {[W]} F.Waldhausen,  {\it  On irreducible 3-manifolds which are
sufficiently large.}  Ann. Math. 87, 56-88  (1968)

\item {[Z1]} B.Zimmermann, {\it On hyperbolic knots with homeomorphic cyclic
branched coverings.} Math. Ann. 311, 665-673 (1998)

\item {[Z2]} B.Zimmermann, {\it On the number and geometry of hyperbolic knots
with a common cyclic branched covering: known results and open problems.}  Rend.
Ist. Mat. Univ. Trieste 38,  95-108  (2006)

\bye